\documentclass[a4paper]{article}

\usepackage{geometry}
\usepackage{pspicture}
\usepackage{pstricks,pst-grad}
\usepackage{amsmath}
\usepackage{amsthm}
\usepackage{amsfonts}
\usepackage{graphicx}
\usepackage{subcaption}
\usepackage{xcolor}
\usepackage{todonotes}
\usepackage{bm}
\usepackage{mathtools}
\usepackage{setspace}
\usepackage{hyperref}
\usepackage{natbib}

\newtheorem{definition}{Definition}[section]
\newtheorem{theorem}{Theorem}[section]
\newtheorem{lemma}{Lemma}[section]
\newcommand{\bs}{\bar{s}}
\allowdisplaybreaks

\begin{document}

  \title{On the Penalty term for the Mixed Discontinuous Galerkin Finite
    Element Method for the Biharmonic Equation}

\author{%
Balaje Kalyanaraman\thanks{Corresponding Author. Email: Balaje.Kalyanaraman@uon.edu.au}\\[2pt]
School of Mathematical and Physical Sciences,\\ University of
Newcastle, Callaghan, NSW 2300, Australia.\\[6pt]
{\sc and}\\
 Danumjaya P.\thanks{Email: danu@goa.bits-pilani.ac.in}\\[2pt]
Birla Institute of Technology and Science, Pilani, K.K.Birla Goa Campus,\\ NH17B, Zuarinagar, South Goa, Goa 403726, India. 
}

\maketitle

  \begin{abstract}
    In this paper, we present a study on the effect of penalty term
    in the mixed Discontinuous Galerkin Finite Element Method for the
    biharmonic equation proposed by \cite{gudi2008mixed}. The proposed
    mixed Discontinuous Galerkin Method showed sub-optimal rates of
    convergence for piecewise quadratic elements and no significant
    convergence rates for piecewise linear elements. We show that by
    choosing the penalty term proportional to $|e_k|^{-1}$ instead of
    $|e_k|^{-3}$, ensures
    an optimal rate of convergence for the approximation, including
    for piecewise linear elements. Finally, we
    present numerical experiments to validate our theoretical
    results.\\[3pt]
    
  \end{abstract}

\section{Introduction}
Discontinuous Galerkin methods are popular finite element techniques
which use discontinuous polynomials to construct approximate
solutions. The local nature of approximation offers flexibility in
using higher order polynomials with non--uniformity in the degree of
approximation and in adaptive methods. For a review of the various
discontinuous Galerkin methods we refer the reader to
\cite{arnold2002unified, cockburn2012discontinuous,
  riviere2008discontinuous}. For fourth order problems, conforming
methods require imposition of $C^1$--continuity across the inter
element boundaries and are computationally expensive. Several finite
element methods have been proposed including the mixed finite element
method \citep{ciarlet1974mixed, monk1987mixed, danumjaya2012mixed},
$C^0$--interior penalty methods \citep{brenner2005c, engel2002continuous} and so on to relax the continuity
requirements. The idea behind the mixed finite element method is to
split the fourth order problem into two second order problems by
introducing an auxiliary variable $v = -\Delta u$. The
$C^1$--continuity requirement is relaxed and the system is then
approximated using a $C^0$-finite element method.\\

An \textit{hp} mixed discontinuous Galerkin Finite Element
Method was proposed by \cite{gudi2008mixed}. Using a
primal formulation leads to integrals involving higher order
derivatives which can be avoided by using a mixed
formulation. However, the method yielded sub-optimal convergence rates
for piecewise quadratic elements and no significant convergence rates
for piecewise linear elements. This is due to the choice of the
penalty parameter in the weak formulation which was taken proportional
to the inverse cube $(|e_k|^{-3})$ of edge/face diameter $e_k$, a common choice in
super-penalization in non-symmetric interior penalty Galerkin
methods \citep{riviere2008discontinuous}. However, modifying the
penalty parameter by taking it proportional to the inverse of the edge
length $(|e_k|^{-1})$, essentially reducing the size of the penalty
term, yields optimal convergence.\\

In this paper, we consider the mixed formulation of the biharmonic
equation
\begin{equation}
  \Delta^2 u = f\quad \text{in}\;\;\Omega,\label{eq:bh}
\end{equation}
subject to the clamped boundary conditions
\begin{equation}
  u = \frac{\partial u}{\partial n} = 0\quad \text{on}\;\;\partial
  \Omega.\nonumber
\end{equation}
where $\Omega \subset \mathbb{R}^n$, $n=2,3$ is a bounded and convex
domain with smooth boundary $\partial \Omega$ and $n$ is the outward
unit normal to $\partial \Omega$. We assume that the data $f$ is
sufficiently smooth so that there exists a unique solution to the
problem in $H^4(\Omega)$. We introduce a new variable $v = -\Delta u$
and split the biharmonic equation~\eqref{eq:bh} into two equations as
\begin{equation}
  \begin{aligned}
    -\Delta v &= f \quad \text{in}\;\;\Omega,\\
    -\Delta u &= v \quad \text{in}\;\;\Omega,\\
    u &= \frac{\partial u}{\partial n} = 0 \quad \text{on} \;\; \Gamma.
  \end{aligned}\label{eq:bh2}
\end{equation}
We then discretize the problem using the mixed discontiunous Galerkin
method and prove some error estimates for the finite element
solution to study the convergence of the method with respect to the
mesh size $h$. We then perform some numerical experiments to validate the
theoretical results.\\

The paper is organized as follows. In Section~\ref{sec:1} we derive
the weak formulation of the biharmonic problem. The proof of
well--posedness of the problem is identical to \cite{gudi2008mixed} and
hence we omit the same. In Section~\ref{sec:2}, we discuss the error
estimates for $h-$refinement and we derive $L^2$ and energy estimates
for the primal variable and an error estimate in the $L^2$-norm for
the auxiliary variable. Finally, in Section~\ref{sec:3} we perform numerical
experiments to validate our theoretical results established in
Section~\ref{sec:2}.

\section{Weak Formulation}
\label{sec:1}
Let $\bar{\Omega} = \underset{K \in \mathcal{T}_h}{\cup}\bar{K}$ be a
uniform partition of $\Omega$. Let us denote the edges of
$\mathcal{T}_h$ by $e_k$ where $k = 1,\cdots,M_h$. Let
$\Gamma_{\text{I}} = \left\{ e_1, e_2, \cdots e_{P_h}\right\}$ denote
the set of all interior edges and $\Gamma_D = \left\{ e_{P_h+1},
  \cdots, e_{M_h} \right\}$, the boundary edges. Let $\Gamma =
\Gamma_{\text{I}} \cup \Gamma_{\text{D}}$. For each edge $e_k = K_i
\cap K_j$, associate a unit normal $n_k$ outward from $K_i$, so that
$n|_{K_i} = -n|_{K_j}$. Let $h_K$ denote the diameter of each $K \in
\mathcal{T}_h$ and $h = \underset{K \in
  \mathcal{T}_h}{\max}\,h_K$. We define the broken Sobolev space
\begin{equation}
  H^s(\mathcal{T}_h) = \left\{ v \in L^2(\Omega) \;:\; v_{K_i}
    \in H^s(K) \;\,\forall K \in \mathcal{T}_h\right\}.
\end{equation}
where $H^s(K)$ denotes the standard Sobolev space of order $s$. The
associated broken norm and semi--norm are defined respectively, by
\begin{equation*}
  \|v\|_{s,\mathcal{T}_h} = \left(\sum_{K \in \mathcal{T}_h}
    \|v\|_{s,K}^2\right)^{1/2}\quad \text{and} \quad
  |v|_{s,\mathcal{T}_h} = \left( \sum_{K \in \mathcal{T}_h}
    |v|_{s,K}^2\right)^{1/2}
\end{equation*}
where $\|v\|_{s,K}$ and $|v|_{s,K}$ denotes the standard Sobolev norm
and seminorm on $K$, respectively. We denote the $L^2$ norm by
$\|\cdot\|$ and its inner--product by $(\cdot,\cdot)$. We set $V =
H^2(\mathcal{T}_h)$.\\

Let $u,v$ be sufficiently smooth functions. We consider the first
equation in \eqref{eq:bh2}, multiply by some $\varphi \in V$, and
integrate over the domain to obtain
\begin{equation}
  -\int_{\Omega} (\Delta v) \varphi \,dx = \int_{\Omega} f\varphi\,dx
\end{equation}
Applying integration by parts and using the fact that $[v] = 0$ in
$\Gamma_{\text{I}}$, we obtain
\begin{equation}
  \begin{aligned}
    \sum_{K \in \mathcal{T}_h} \int_K \nabla v\cdot \nabla \varphi\,dx -
    \sum_{e_k \in \Gamma}\int_{e_k} \left\{ \nabla v\cdot
      n_k\right\}[\varphi]\,ds &- \sum_{e_k \in
      \Gamma_{\text{I}}}\int_{e_k} \left\{ \nabla \varphi\cdot
      n_k\right\}[v]\,ds\\
    &= \sum_{K \in\mathcal{T}_h}\int_K f\varphi\, dx.\label{eq:3}
  \end{aligned}
\end{equation}
Since $[u] = 0$ in $\Gamma$, we see that
\begin{equation}
  J(u,\varphi) = \sum_{e_k \in \Gamma}\int_{e_k}
  \alpha_k[u][\varphi]\, ds = 0 \label{eq:4}
\end{equation}
where $\alpha_k$ is a positive real constant. Adding
$J(u,\varphi)$ to~\eqref{eq:3} we obtain
\begin{equation}
  \begin{aligned}
    \sum_{K \in \mathcal{T}_h} \int_K \nabla v\cdot \nabla \varphi\,dx &-
    \sum_{e_k \in \Gamma}\int_{e_k} \left\{ \nabla v\cdot
      n_k\right\}[\varphi]\,ds - \sum_{e_k \in
      \Gamma_{\text{I}}}\int_{e_k} \left\{ \nabla \varphi\cdot
      n_k\right\}[v]\,ds\\
    &+ \sum_{e_k \in \Gamma}\int_{e_k}
    \alpha_k[u][\varphi]\, ds
    = \sum_{K \in\mathcal{T}_h}\int_K f\varphi\, dx,\label{eq:5}
  \end{aligned}
\end{equation}
which yields the weak formulation of the first equation. Similarly for
the second equation, multiplying some $\chi \in V$, integrating over
$\Omega$ and using the fact that $\nabla u\cdot n = 0$ on $\Gamma_{\text{D}}$,
$[u] = 0$ in $\Gamma$, we obtain
\begin{equation}
  \begin{aligned}
    \sum_{K \in \mathcal{T}_h} \int_K \nabla u\cdot \nabla \chi\,dx &-
    \sum_{e_k \in \Gamma_{\text{I}}}\int_{e_k} \left\{ \nabla u\cdot
      n_k\right\}[\chi]\,ds - \sum_{e_k \in
      \Gamma}\int_{e_k} \left\{ \nabla \chi\cdot
      n_k\right\}[u]\,ds
    = \sum_{K \in\mathcal{T}_h}\int_K v\chi\, dx,\label{eq:6}
  \end{aligned}
\end{equation}
We define the bilinear forms
\begin{equation*}
  \begin{aligned}
    B(w,z) &= \sum_{K\in \mathcal{T}_h}\int_K \nabla w\cdot \nabla
    z\,dx - \sum_{e_k \in \Gamma}\int_{e_k}\left\{ \nabla
      w\cdot n_k\right\}[z]\,ds - \sum_{e_k \in
      \Gamma_{\text{I}}}\int_{e_k}\left\{ \nabla
      z\cdot n_k\right\}[w]\,ds\\
    B_{\text{I}}(w,z) &= \sum_{K\in \mathcal{T}_h}\int_K \nabla w\cdot \nabla
    z\,dx - \sum_{e_k \in \Gamma_{\text{I}}}\int_{e_k}\left\{ \nabla
      w\cdot n_k\right\}[z]\,ds - \sum_{e_k \in
      \Gamma_{\text{I}}}\int_{e_k}\left\{ \nabla
      z\cdot n_k\right\}[w]\,ds\\
  \end{aligned}
\end{equation*}
Now the weak formulation of the problem is as follows: Find $(u,v) \in
V \times V$, such that
\begin{equation}
  \begin{aligned}
    B(v,\varphi) + J(u,\varphi) &= (f,\varphi),\\
    B(\chi, u) &= (v,\chi),
  \end{aligned}\label{eq:7}
\end{equation}
for all $(\varphi,\chi) \in V \times V$. Define a finite dimensional subspace
$V_h$ of $V$ as
\begin{equation*}
  V_h = \left\{ v \in H^s(\mathcal{T}_h) \; : \; v|_K
    \in \mathbb{P}_p(K)\;\; \text{for all}\;\; K \in \mathcal{T}_h \right\},
\end{equation*}
where $\mathbb{P}_p(K)$ denotes the space of polynomials of degree
$\le p$. Now the discontinuous Galerkin finite element weak
formulation reads: Find $(u_h,v_h) \in V_h \times V_h$ such that
\begin{equation}
  \begin{aligned}
    B(v_h,\varphi) + J(u_h,\varphi) &= (f,\varphi),\\
    B(\chi, u_h) &= (v_h,\chi),
  \end{aligned}\label{eq:8}
\end{equation}
for all $(\varphi,\chi) \in V_h \times V_h$. For the purpose of error
analysis, we define the following mesh dependent energy norms,
\begin{equation*}
  \begin{aligned}
    |||w|||^2 &= \left( \sum_{K\in\mathcal{T}_h} \int_K |\nabla w|^2
      \,dx + \sum_{e_k \in \Gamma} \int_{e_k} \frac{|e_k|}{\sigma_1}
      \left\{\frac{\partial w}{\partial n_k}\right\}^2 \, ds+
      \sum_{e_k\in \Gamma}\int_{e_k} \frac{\sigma_1}{|e_k|}
      [w]^2\,ds\right),\\
    |||w|||_{\text{I}}^2 &= \left( \sum_{K\in\mathcal{T}_h} \int_K |\nabla w|^2
      \,dx + \sum_{e_k \in \Gamma} \int_{e_k} \frac{|e_k|}{\sigma_1}
      \left\{\frac{\partial w}{\partial n_k}\right\}^2 \, ds+
      \sum_{e_k\in \Gamma_{\text{I}}}\int_{e_k} \frac{\sigma_1}{|e_k|} [w]^2\,ds\right).
  \end{aligned}
\end{equation*}
Below, we state some properties of the bilinear form
$B(\cdot,\cdot)$ without proof.
\begin{lemma}
  \label{lemma1}
  For sufficiently large constant $\sigma_1$, it can
  be shown that for any $w_h \in V_h$,
  \begin{equation}
    C|||w_h|||_{\text{I}}^2 \le B_{\text{I}}(w_h,w_h) + \sum_{e_k \in \Gamma_{\text{I}}}
    \int_{e_k} \frac{\sigma_1}{|e_k|} [w_h]^2\,ds,\label{eq:coercivity}
  \end{equation}
  where $C$ is a constant independent of $h$ and $p$. It can also be
  shown that for all $w,q \in V$ there exists a constant independent of
  $h$ such that
  \begin{equation}
    |B(w,q)| \le C|||w|||_{\text{I}}\;|||q|||.\label{eq:bounds}
  \end{equation}
\end{lemma}
The inequalities \eqref{eq:coercivity}~and~\eqref{eq:bounds} refer to
the coercivity and the boundedness property of the bilinear forms. For
proofs, we refer the reader to
\cite{riviere2008discontinuous,gudi2008mixed} and the references
therein. Below we define an interpolation operator for functions on
the Broken Sobolev space.
\begin{definition}
  For given $\phi \in H^s(\mathcal{T}_h)$, define $I_h\phi \in V_h$
  be an interpolation operator satisfying optimal approximation
  properties, i.e.,
  \begin{equation*}
    |||\phi - I_h\phi||| \le Ch^{\mu-1}\|\phi\|_{s,\mathcal{T}_h}
  \end{equation*}
  where $\mu=\min\{p+1,s\}$.
\end{definition}
We recall the following trace inequality on the finite
element space.
\begin{lemma}
  \label{lemma2}
  Let $v_h \in V_h$. There exists a constant $C > 0$ such that
  \begin{equation}
    \|\nabla^l v_h\|_{L^2(e_k)} \le Ch_K^{-1/2}\|\nabla^l
    v_h\|_{L^2(K)}, \quad l=0,1.\label{eq:trace}
  \end{equation}
\end{lemma}
Below, we state the inverse inequality without proof.
\begin{lemma}
  \label{lemma3}
  Let $v_h \in V_h$. There exists a constant $C > 0$ such that
  \begin{equation}
    |v_h|_{H^1(K)} \le Ch_K^{-1}\|v_h\|_{L^2(K)},\label{eq:inverse}
  \end{equation}
\end{lemma}
For proofs of the trace
inequality and the inverse inequality, we refer the reader to
\cite{brenner2007mathematical} and
\cite{riviere2008discontinuous}. 
We state and prove the following Lemma.
\begin{lemma}
  \label{lemma4.1}
  For $w_h \in V_h$, there exists a positive constant $C>0$, such that
  \begin{equation*}
    |||w_h|||^2 \le C_T\left(\sum_{K\in \mathcal{T}_h} \int_K |\nabla
      w_h|^2 + \sum_{e_k \in \Gamma}\int_{e_k} \frac{\sigma_1}{|e_k|}[w_h]^2\,ds\right).
  \end{equation*}
  \begin{proof}
    From the definition of the energy norm and the trace inequality,
    we have for $w_h \in V_h$
    \begin{align*}
      |||w_h|||^2 &= \left( \sum_{K\in\mathcal{T}_h} \int_K |\nabla w_h|^2
                    \,dx + \sum_{e_k \in \Gamma} \int_{e_k} \frac{|e_k|}{\sigma_1}
                    \left\{\frac{\partial w_h}{\partial n_k}\right\}^2 \, ds+
                    \sum_{e_k\in \Gamma}\int_{e_k} \frac{\sigma_1}{|e_k|}
                    [w_h]^2\,ds\right),\\
                  &\le C_T\left(\sum_{K\in\mathcal{T}_h} \int_K |\nabla w_h|^2
                    \,dx + \sum_{e_k\in \Gamma}\int_{e_k}
                    \frac{\sigma_1}{|e_k|}
                    [w_h]^2\,ds\right),
    \end{align*}
    which yields the desired result.
  \end{proof}
\end{lemma}

We are interested to study the effect of the penalty parameter on the
convergence of the discrete solution. Let us set
\begin{equation}
  \alpha_{k} = \sigma_0 |e_k|^{-i} p^2. \label{eq:penalty}
\end{equation}
where $p$ is the degree of approximation of the polynomial,
$\sigma_0$ is a constant which will be defined later and $i$ is
an integer which describes the degree of penalization. It
is well known in literature that for the case of the Non--symmetric
interior penalty Galerkin (NIPG) method for second order elliptic
problems, the choice of the penalty parameter plays a crucial role on
the convergence of the solution. If $i=1$, under normal penalization,
the convergence in the $L^2-$norm for the piecewise quadratic case is
suboptimal. However, for the case when $i=3$, under super
penalization, the convergence is optimal.
For the mixed Discontinuous Galerkin Method considered here, we observe that the
penalization term with $i=1$ produces optimal convergence in all
cases, as opposed to $i=3$ where optimal convergence rates are
observed only for piecewise cubic elements. The latter case was well
studied by \cite{gudi2008mixed}. In the next section, we
derive optimal error estimates with $i=1$.

\section{Error Analysis}
\label{sec:2}
Define the auxiliary projection $\Pi_h : H^s(\mathcal{T}_h) \to
V_h$ by
\begin{equation}
  \begin{aligned}
    B_{\text{I}}(\phi - \Pi_h \phi, \chi) &+ \sum_{e_k \in \Gamma_{\text{I}}}
    \int_{e_k} \frac{\sigma_1}{|e_k|}[\phi - \Pi_h \phi][\chi]\,ds =
    0,\label{eq:Proj}\\
    \left(\phi - \Pi_h \phi, 1\right) &= 0.
  \end{aligned}
\end{equation}
The auxiliary projection can be shown to satisfy optimal error
estimates in the $L^2$ and energy norms. We define
$\mu=\min\{p+1,s\}$, $\bar{\mu}=\min\{p+1,\bs\}$,
$\theta=\min\{p+1,4\}$. We state and prove the
following Lemma.
\begin{lemma}
  \label{lemma4}
  For a given $\phi \in V$, there exists a unique $\Pi_h\phi \in V_h$
  satisfying \eqref{eq:Proj}. Moreover, there exists positive
  constants independent of $h$ such that
  \begin{align}
    |||\phi - \Pi_h\phi||| &\le Ch^{\mu-1}\,\|\phi\|_{s,\mathcal{T}_h}\\
    \|\phi - \Pi_h\phi\| &\le Ch^{\mu}\,\|\phi\|_{s,\mathcal{T}_h}.
  \end{align}
  provided $\sigma_1$ is a sufficiently large positive constant, i.e., $\sigma_1 \ge \sigma_1^*$ for some $\sigma_1^* \in \mathbb{R}^+$.
  \begin{proof}
    For the proof, we refer the reader to Theorem 5.1 in \cite{gudi2008mixed}.
  \end{proof}
\end{lemma}
Subtracting \eqref{eq:8} from
\eqref{eq:7}, we obtain the error equations
\begin{align}
  B(v - v_h, \varphi) + J(u - u_h, \varphi) &= 0,\label{eq:orth1}\\
  B(\chi, u - u_h) - (v - v_h,\chi) &= 0.\label{eq:orth2}
\end{align}
for $(\varphi,\chi) \in V_h \times V_h$. Set $e_u:= u-u_h$ and $e_v :=
v-v_h$. The following Lemma is useful in proving the error estimates.
\begin{lemma}
  \label{lemma3.1}
  Let $\xi_v = v_h - \Pi_h v$. Provided $\sigma_0 > \sigma_*$, there
  exists a positive constant
  $C_*$ such that
  \begin{align*}
    C_*|||\xi_v|||^2 \le Ch^{2\mu-2} \|u\|_{s,\mathcal{T}_h}^2 +
    Ch^{2\bar{\mu}-2} \|u\|_{\bs+2,\mathcal{T}_h}^2 + \frac{\sigma_0
    \epsilon^2}{2C_1} J(\xi_u,\xi_u)
  \end{align*}
  where $C_1$ is a positive constant independent of $\sigma_0$ and $0
  < \epsilon \ll 1$.
  \begin{proof}
    Let $\eta_v = v - \Pi_h v$. From equation~\eqref{eq:orth1} we obtain,
    \begin{equation*}
      B(\xi_v,\varphi) = B(\eta_v,\varphi) + J(\eta_u,\varphi) - J(\xi_u,\varphi),
    \end{equation*}
    where $\eta_u = u - I_h u$ and $\xi_u = u_h - I_h u$. Setting
    $\varphi = \xi_v$ and expanding the bilinear form
    $B(\cdot,\cdot)$, we obtain
    \begin{align}
      \sum_{K\in \mathcal{T}_h} \int_K |\nabla \xi_v|^2 \,dx
      &=
        \sum_{e_k \in \Gamma_\text{I}}\int_{e_k} \left\{\frac{\partial
        \xi_v}{\partial n_k}\right\}[\xi_v]\,ds + \sum_{e_k \in \Gamma}\int_{e_k}
        \left\{\frac{\partial \xi_v}{\partial
        n_k}\right\}[\xi_v]\,ds\nonumber\\
      &+
        B(\eta_v,\xi_v) + J(\eta_u,\xi_v) - J(\xi_u,\xi_v).\label{eq:newErr1}
    \end{align}
    Using Lemma~\ref{lemma4.1} with $w_h =\xi_v$, and substituting the
    result of equation~\eqref{eq:newErr1} onto the Lemma, we obtain the inequality
    \begin{align*}
      c|||\xi_v|||^2
      &\le \sum_{e_k \in \Gamma_\text{I}}\int_{e_k} \left\{\frac{\partial
        \xi_v}{\partial n_k}\right\}[\xi_v]\,ds + \sum_{e_k \in \Gamma}\int_{e_k}
        \left\{\frac{\partial \xi_v}{\partial
        n_k}\right\}[\xi_v]\,ds\nonumber\\
      &+
        B(\eta_v,\xi_v) + J(\eta_u,\xi_v) - J(\xi_u,\xi_v) + \sum_{e_k
        \in \Gamma} \int_{e_k} \frac{\sigma_1}{|e_k|} [\xi_v]^2\,ds,
    \end{align*}
    with $c=\frac{1}{C_T}$ being a constant which arises due to the trace inequality and is dependent on the triangulation. We observe that using the boundedness properties of the bilinear forms
    together with the trace and Young's inequalities, we have
    \begin{align*}
      B(\eta_v,\xi_v) &\le C|||\eta_v|||\;|||\xi_v||| \le
                        Ch^{2\bar{\mu}-2}\|u\|_{\bs+2,\mathcal{T}_h}^2
                        + \frac{\epsilon^2}{4}|||\xi_v|||^2,\\
      J(\eta_u,\xi_v) &\le (J(\eta_u,\eta_u))^{1/2}
                        (J(\xi_v,\xi_v))^{1/2} \le Ch^{2\mu-2}
                        \|u\|_{s,\mathcal{T}_h}^2 +
                        \frac{\epsilon^2}{4}|||\xi_v|||^2,
    \end{align*}
    Using the Cauchy-Schwarz's inequality and
    subsequently the Young's inequality, the trace inequality and the inverse
    inequality, the term,
    \begin{align*}
      \sum_{e_k \in \Gamma}\int_{e_k}
      \left\{\frac{\partial \xi_v}{\partial
      n_k}\right\}[\xi_v]\,ds &\le
                                \left(\sum_{e_k\in
                                \Gamma}\int_{e_k}\frac{|e_k|}{\sigma_0}\left\{\frac{\partial
                                \xi_v}{\partial
                                n_k}\right\}\,ds\right)^{1/2}
                                \left(\sum_{e_k\in \Gamma}\int_{e_k}
                                \frac{\sigma_0}{|e_k|}
                                [\xi_v]^2\,ds\right)^{1/2},\\
                              &\le \left(\sum_{K\in\mathcal{T}_h}
                                \int_K \frac{1}{\sigma_0}|\nabla
                                \xi_v|^2\,dx\right)^{1/2}
                                \left(\sum_{e_k\in \Gamma}\int_{e_k}
                                \frac{\sigma_0}{|e_k|}
                                [\xi_v]^2\,ds\right)^{1/2},\\
                                &\le \frac{C_1}{\sigma_0} |||\xi_v||| \; |||\xi_v|||,\\
                              &\le \frac{C_1}{2\sigma_0\epsilon^2}|||\xi_v|||^2 +
                                \frac{\epsilon^2}{2}|||\xi_v|||^2,
    \end{align*}
    where $C_1$ is a constant independent of $\sigma_0$ and depends on the triangulation. Similarly, we have
    \begin{equation*}
        J(\xi_u,\xi_v) \le \frac{\sigma_0 \epsilon^2}{2C_1} J(\xi_u,\xi_u)
                       + \frac{C_1}{2\sigma_0\epsilon^2}|||\xi_v|||^2.
    \end{equation*}
    Finally we have the bound,
    \begin{equation*}
      \sum_{e_k
        \in \Gamma} \int_{e_k} \frac{\sigma_1}{|e_k|} [\xi_v]^2\,ds
      \le \frac{\sigma_1}{\sigma_0} |||\xi_v|||^2.
    \end{equation*}
    Combining all the bounds for a small $\epsilon$, we obtain the intermediate inequality,
    \begin{align*}
      \left(c - \frac{C_1}{\sigma_0\epsilon^2} -
      \frac{\sigma_1}{\sigma_0}\right)|||\xi_v|||^2 \le
      Ch^{2\mu-2}\|u\|_{s,\mathcal{T}_h}^2 +
      Ch^{2\bar{\mu}-2}\|u\|_{\bs+2,\mathcal{T}_h} + \frac{\sigma_0\epsilon^2}{2C_1}J(\xi_u,\xi_u)
    \end{align*}
    A condition that
    \begin{align*}
      C_*=\left(c - \frac{C_1}{\sigma_0\epsilon^2} -
      \frac{\sigma_1}{\sigma_0}\right) > 0 \quad \text{or} \quad
      \sigma_0 > \frac{C_1+\sigma_1\epsilon^2}{c\epsilon^2},
    \end{align*}
    yields a lower bound on $\sigma_0$. Thus,
    \begin{align*}
      C_*|||\xi_v|||^2 \le
      Ch^{2\mu-2}\|u\|_{s,\mathcal{T}_h}^2 +
      Ch^{2\bar{\mu}-2}\|u\|_{\bs+2,\mathcal{T}_h} + \frac{\sigma_0\epsilon^2}{2C_1}J(\xi_u,\xi_u),
    \end{align*}
    which is the desired estimate with $\sigma_* =
    \frac{C_1+\sigma_1\epsilon^2}{c\epsilon^2}$ for a sufficiently small $\epsilon$.
  \end{proof}
\end{lemma}
Lemma~\ref{lemma3.1} relates the constant $\sigma_0$ with $\sigma_1$
by providing a bound on $\sigma_0$ in terms of $\sigma_1$ which in
turn depends on the triangulation and is bounded below by
$\sigma_1^*$, which is independently known from $\sigma_0$. Now we
state and prove the following
theorems which yields the error
estimates for the DG method.
\begin{theorem}
  \label{theorem1}
  Let $u_h,v_h$ satisfy \eqref{eq:8}. Then there exists some positive
  constant $C$ such that
  \begin{align*}
    \|e_v\|^2 + J(e_u,e_u) \le C\left( h^{2\mu-2}
    \|u\|_{s,\mathcal{T}} +
    h^{2\bar{\mu}-2}\|u\|_{\bar{s},\mathcal{T}}\right)
  \end{align*}
  provided $\sigma_* < \sigma_0 < \sigma^*$ for some
  $\sigma_*,\sigma^* \in \mathbb{R}^+$.
  \begin{proof}
    Let $\delta \ll 1$. First, we split
    \begin{align*}
      u - u_h &= (u - I_h u) - (u_h - I_h u) :=
                \eta_u - \xi_u,\\
      v - v_h &= (v - \Pi_h v) - (v_h - \Pi_h v) :=
                \eta_v - \xi_v.
    \end{align*}
    Substituting the above definitions into the error equations
    \eqref{eq:orth1}~and~\eqref{eq:orth2} setting $\varphi=\xi_u$ in
    \eqref{eq:orth1} and $\chi=\xi_v$ in \eqref{eq:orth2} and subtracting the resulting equations, we obtain
    \begin{align*}
      \|\xi_v\|^2 + J(\xi_u,\xi_u) = B(\eta_v,\xi_u) +
      J(\eta_u,\xi_u) - B(\xi_v,\eta_u) + (\eta_v,\xi_v).
    \end{align*}
    Using the definition of the auxiliary projection we obtain,
    \begin{align*}
      B(\eta_v,\xi_u) &= B_{\text{I}}(\eta_v,\xi_u) - \sum_{e_k\in
                        \Gamma_{\text{D}}} \int_{e_k} \frac{\partial \eta_v}{\partial
                        n_k}\xi_u\,ds,\\
                      &= -\sum_{e_k \in \Gamma_{\text{I}}} \int
                        \frac{\sigma_1}{|e_k|} [\eta_v][\xi_u]\,ds- \sum_{e_k\in
                        \Gamma_{\text{D}}} \int_{e_k} \frac{\partial \eta_v}{\partial
                        n_k}\xi_u\,ds.
    \end{align*}
    We observe that using the trace inequality, inverse inequality and
    the Young's inequality,
    \begin{align*}
      B(\eta_v,\xi_u) \le Ch^{2\bar{\mu}-2}\|u\|_{\bs+2,\mathcal{T}_h}^2 +
      \frac{\delta^2}{2} J(\xi_u,\xi_u).
    \end{align*}
    for some $\delta \ll 1$. Using the Cauchy-Schwarz inequality and the Young's inequality, we obtain
    \begin{align*}
      J(\eta_u,\xi_u) \le Ch^{2\mu-2}\|u\|_{s,\mathcal{T}_h}^2 + \frac{\delta^2}{2} J(\xi_u,\xi_u).
    \end{align*}
    We observe that using the Young's inequality and the result from
    Lemma~\ref{lemma3.1} for the estimate of $|||\xi_v|||$ we obtain,
    \begin{align*}
      B(\xi_v,\eta_u) &\le C|||\xi_v|||\;|||\eta_u|||,\\
                      &\le Ch^{2\mu-2}\|u\|_{s,\mathcal{T}_h}^2 +
                        \frac{\delta^2 C_*}{2} |||\xi_v|||^2,\\
                      &\le Ch^{2\mu-2}\|u\|_{s,\mathcal{T}_h}^2 +
                        Ch^{2\bar{\mu}-2} \|u\|_{\bs+2,\mathcal{T}_h}^2 +
                        \frac{\sigma_0\epsilon^2\delta^2}{4C_1}J(\xi_u,\xi_u),
    \end{align*}
    provided $\sigma_0 > \sigma_*$. We note that the constant $C_*$
    which appears in front of $|||\eta_v|||$ is absorbed into the
    estimate involving $h$. Since $C_*$
    does not depend on $h$, the estimate is unaffected. Finally, we observe that
    \begin{equation*}
      (\eta_v,\xi_v) \le Ch^{2\bar{\mu}}\|u\|_{\bs+2,\mathcal{T}_h} + \frac{\delta^2}{2}\|\xi_v\|^2.
    \end{equation*}
    Combining the bounds and
    using the triangle inequality, we obtain
    \begin{equation}
      \|\xi_v\|^2 + J(\xi_u,\xi_u) \le Ch^{2\mu-2}
      \|u\|_{s,\mathcal{T}_h} +
      Ch^{2\bar{\mu}-2}\|u\|_{\bar{s}+2,\mathcal{T}_h} +
      \frac{\sigma_0\epsilon^2\delta^2}{4C_1}J(\xi_u,\xi_u).
    \end{equation}
    Setting
    \begin{align*}
      1 -\frac{\sigma_0\epsilon^2\delta^2}{4C_1} > 0 \quad
      \text{or} \quad \sigma_0 < \frac{4C_1}{\epsilon^2\delta^2},
    \end{align*}
    yields an upper bound for $\sigma_0$ in terms of the constant
    $C_1$. The desired estimate follows from the triangle inequality and with $\sigma_*=
    \frac{C_1+\sigma_1\epsilon^2}{c\epsilon^2}$ and
    $\sigma^*=\frac{4C_1}{\epsilon^2\delta^2}$ for sufficiently small $\delta$. 
  \end{proof}
\end{theorem}
The consequence of having the condition that the constant $\sigma_0$ to be bounded between
$[\sigma_*,\sigma^*]$ may explain why the convergence rate
deteriorates for large values of penalty parameters. Choosing the value of $\sigma_0$ too large tends the
value of the penalty terms to be proportional to $|e_k|^{-3}$ and thusumay erode the estimates derived above. With the above conditions, the estimate in Theorem~\ref{theorem1} predicts a higher rate of convergence, $\mu-1$ for $\|e_v\|$ than $\mu-2$ which was predicted by \cite{gudi2008mixed}. This will lead to better estimates for $|||e_u|||$ and $\|e_v\|$. We state and prove the following theorems.
\begin{theorem}
  \label{theorem2}
  Let $u_h,v_h$ satisfy \eqref{eq:8}. Then there exists some positive
  constant $C$ such that
  \begin{align*}
    \|e_u\| \le
    C\left( h|||e_u||| + h^2\|e_v\| +
    h^{\theta-1}|||v-\Pi_h v||| + h^{\theta-2}\|\Pi_h v-v_h\| +
    h^{\theta-1}(J(e_u,e_u))^{1/2}\right).
  \end{align*}
  \begin{proof}
    For the construction of the proof, we use the Aubin-Nitsche
    duality argument. Consider the dual problem,
    \begin{align*}
      -\Delta \phi &= z \quad \text{in} \quad \Omega,\\
      -\Delta \psi &= \phi \quad \text{in} \quad \Omega,\\
      \psi &= \frac{\partial \psi}{\partial n} = 0 \quad \text{on} \quad \partial\Omega,
    \end{align*}
    where the functions $\phi$ and $\psi$ satisfy the regularity result
    \begin{align*}
      \|\psi\|_{H^4} + \|\phi\|_{H^2} \le C\|z\|.
    \end{align*}
    Take $z=e_u$, multiply the first equation of the dual problem by
    $e_u$ and integrate over the domain to obtain
    \begin{equation}
      B(\phi,e_u) + J(e_u,\psi) = \|e_u\|^2.\label{eq:error3}
    \end{equation}
    where we have used the fact that $[\psi]=0$ on $\Gamma$ and $[\phi]=0$ on $\Gamma_{\text{I}}$. Next multiply the second equation of the dual problem by $e_v$ and
    integrate to obtain
    \begin{equation}
      B(e_v,\psi) - (e_v,\phi) = 0.\label{eq:error4}
    \end{equation}
    where we have used $\frac{\partial \psi}{\partial n}=0$ on $\Gamma_{\text{D}}$ and $[\psi]=0$ on $\Gamma$.  Adding the equations \eqref{eq:error3} and \eqref{eq:error4} we obtain
    \begin{equation}
      \|e_u\|^2 = B(\phi,e_u) + J(e_u,\psi) + B(e_v,\psi) - (e_v,\phi).\label{eqnewErr}
    \end{equation}
    Denoting $\eta_\phi = \phi - I_h\phi$ and $\eta_\psi = \psi -
    I_h\psi$, from the orthogonality result \eqref{eq:orth1} and
    \eqref{eq:orth2} we have 
    \begin{align*}
        B(e_v, I_h\psi) + J(e_u,I_h\psi) &= 0,\\
        B(I_h\phi,e_u) - (e_v,I_h\phi) &= 0.
    \end{align*}
Subtracting the above equations from \eqref{eqnewErr}, we obtain the following equation
    \begin{equation}
      \|e_u\|^2 = B(\eta_\phi,e_u) + J(e_u,\eta_\psi) +
      B(e_v,\eta_\psi) - (e_v,\eta_\phi).\label{eq:error5}
    \end{equation}
    Constructing the error bounds and using the regularity of the functions
    $\phi$ and $\psi$ we obtain,
    \begin{align*}
      B(\eta_\phi,e_u) &\le C|||\eta_\phi|||\,|||e_u||| \le
                         Ch|||e_u|||\;\|\phi\|_{H^2},\\
      (e_v,\eta_\phi) &\le Ch^2\|e_v\|\;\|\phi\|_{H^2},\\
      B(e_v,\eta_\psi) &\le C|||e_v|||\;|||\eta_\psi||| \le
                         Ch^{\theta-1}|||e_v|||\;\|\psi\|_{H^4},\\
      J(e_u,\eta_\psi) &\le (J(\eta_\psi,\eta_\psi))^{1/2} \,
                         (J(e_u,e_u))^{1/2} \le Ch^{\theta-1}\;(J(e_u,e_u))^{1/2}\|\psi\|_{H^4}.
    \end{align*}
    Combining all the bounds,
    \begin{align}
      \|e_u\| &\le C\left( h|||e_u||| + h^2\|e_v\| +
                h^{\theta-1}|||e_v||| + h^{\theta-1}
                (J(e_u,e_u))^{1/2}\right),\nonumber\\
              &\le C\left( h|||e_u||| + h^2\|e_v\| +
                h^{\theta-1}|||v-\Pi_h v||| + h^{\theta-2}\|\Pi_h v-v_h\| +
                h^{\theta-1}(J(e_u,e_u))^{1/2}\right),\label{eq:error6}
    \end{align}
    which yields the desired estimate.
  \end{proof}
\end{theorem}

\begin{theorem}
  \label{theorem3}
  Let $u_h,v_h$ satisfy \eqref{eq:8}. There exists some positive
  constant $C$ such that
  \begin{align*}
    |||e_u||| &\le
                Ch^{\min\{\mu+\theta-3,\mu-1\}}\|u\|_{s,\mathcal{T}_h}
                + Ch^{\min\{\bar{\mu}+\theta-3,\bar{\mu}-1\}}\|u\|_{\bs+2,\mathcal{T}_h},
  \end{align*}
  provided the conditions on $\sigma_0$ in Theorem~\ref{theorem1} holds.
  \begin{proof}
    Let the constant $\epsilon\ll
    1$ and $\xi_u = u_h - I_h u$. We split the error $u -
    u_h = (u - I_h u) - (u_h - I_h u) := \eta_u - \xi_u$. From Lemma~\ref{lemma4.1}, we note that
    \begin{align*}
      c|||\xi_u|||^2 &\le \left(\sum_{K\in \mathcal{T}_h} \int_K
                       |\nabla \xi_u|^2\,dx + \sum_{e_k\in \Gamma}
                       \int_{e_k}\frac{\sigma_1}{|e_k|}
                       [\xi_u]^2\,ds\right)\\
                     &\le B(\xi_u,\xi_u) + \sum_{e_k\in \Gamma}\int_{e_k}
                       \left\{\frac{\partial \xi_u}{\partial n_k}\right\}[\xi_u]\,ds
                       + \sum_{e_k\in \Gamma_{\text{I}}}\int_{e_k}
                       \left\{\frac{\partial \xi_u}{\partial
                       n_k}\right\}[\xi_u]\,ds +
                       \sum_{e_k\in \Gamma}
                       \int_{e_k}\frac{\sigma_1}{|e_k|}
                       [\xi_u]^2\,ds
    \end{align*}
    From \eqref{eq:orth2}, we obtain
    \begin{equation*}
      B(\chi, \xi_u) = B(\chi,\eta_u) - (e_v,\chi).
    \end{equation*}
    Thus, with $\chi=\xi_u$ we obtain the inequality
    \begin{align*}
      c|||\xi_u|||^2  \le B(\xi_u,\eta_u) - (e_v,\xi_u) &+ \sum_{e_k\in \Gamma}\int_{e_k}
                                                          \left\{\frac{\partial \xi_u}{\partial n_k}\right\}[\xi_u]\,ds
                                                          + \sum_{e_k\in \Gamma_{\text{I}}}\int_{e_k}
                                                          \left\{\frac{\partial \xi_u}{\partial
                                                          n_k}\right\}[\xi_u]\,ds \\
                                                        &+ \sum_{e_k\in \Gamma}
                                                          \int_{e_k}\frac{\sigma_1}{|e_k|}
                                                          [\xi_u]^2\,ds.
    \end{align*}
    We construct the following upper bounds for the terms on the right
    hand side,
    \begin{align*}
      B(\xi_u,\eta_u) &\le C|||\xi_u|||\; |||\eta_u|||| \le
                        Ch^{2\mu-2}\|u\|_{s,\mathcal{T}_h}^2 + \frac{\epsilon^2}{2}
                        |||\xi_u|||^2,\\
      (e_v,\xi_u) &\le \frac{C}{2\epsilon^2}\|e_v\|^2 +
                    \frac{\epsilon^2}{2}\|\xi_u\|^2,\\
      \sum_{e_k \in \Gamma}\int_{e_k} \left\{\frac{\partial
      \xi_u}{\partial n_k}\right\}[\xi_u]\,ds &\le \left( \sum_{e_k \in
                                                \Gamma}\int_{e_k}\frac{|e_k|}{\sigma_0}
                                                \left\{\frac{\partial
                                                \xi_u}{\partial n_k}\right\}^2\,ds\right)^{1/2} \,\left(\sum_{e_k \in
                                                \Gamma}\int_{e_k}\frac{\sigma_0}{|e_k|}
                                                [\xi_u]^2\,ds\right)^{1/2},\\
                      &\le \frac{C}{2\epsilon^2}J(\xi_u,\xi_u) + \frac{\epsilon^2}{2} |||\xi_u|||^2.
    \end{align*}
    Fixing $\sigma_0$ and noting that $k_1\sigma_0 \le \sigma_1 \le k_2 \sigma_0$, we have
    the final upper bound,
    \begin{equation*}
      \sum_{e_k\in \Gamma} \int_{e_k}\frac{\sigma_1}{|e_k|}
      [\xi_u]^2\,ds \le k_2 J(\xi_u,\xi_u).
    \end{equation*}
    Combining all the bounds and applying the triangle inequality, we obtain
    \begin{align}
      c|||\xi_u|||^2 \le Ch^{2\mu-2}
      \|u\|_{s,\mathcal{T}_h}^2 + C\left(\|e_v\|^2 + J(\xi_u,\xi_u)
      \right) + \frac{\epsilon^2}{2}\|\xi_u\|^2. \label{eq:error7}
    \end{align}

    \noindent
    Writing
    \begin{align*}
      e_v = v - v_h = v - \Pi_h v - (v_h - \Pi_h v),
    \end{align*}
    we obtain $\Pi_h v - v_h = (v - \Pi_h v) - e_v$. From estimate~\eqref{eq:error6}, we have
    \begin{align*}
      \epsilon \|e_u\| &\le C\epsilon\left( h|||e_u||| + h^2\|e_v\| +
                         h^{\theta-1}|||\Pi_h v - v||| + h^{\theta-2}\|\Pi_h v-v_h\| +
                         h^{\theta-1}(J(e_u,e_u))^{1/2}\right),\\
                       &\le C\epsilon\left(h|||e_u|||
                         + h^2\|e_v\| + h^{\theta-1}|||v-\Pi_hv|||+ h^{\theta-2} \|v-\Pi_h v\| +
                         h^{\theta-2}\left(\|e_v\| + J(e_u,e_u)^{1/2}\right)\right).
    \end{align*}
    Now we have from the results of Theorem~\ref{theorem1} and
    Lemma~\ref{lemma4} the estimates,
    \begin{align*}
      h^2 \|e_v\| + h^{\theta-2}\left(\|e_v\| +
      J(e_u,e_u)^{1/2}\right) &\le
                                Ch^{\min\{\theta-2,2\}}\left(h^{\mu-1} \|u\|_{s,\mathcal{T}_h} +
                                h^{\bar{\mu}-1}
                                \|u\|_{\bs+2,\mathcal{T}_h} \right),\\
      h^{\theta-2}\|v-\Pi_h v\| &\le Ch^{\bar{\mu}+\theta-2}\|u\|_{\bs+2,\mathcal{T}_h},\\
      h^{\theta-1}|||v-\Pi_h v||| &\le Ch^{\bar{\mu}+\theta-2}\|u\|_{\bs+2,\mathcal{T}_h},
    \end{align*}
    respectively. Thus the estimate for $e_u$ reads,
    \begin{align}
      \epsilon\|e_u\| \le C\epsilon \left( h|||e_u||| +
      h^{\bar{\mu}+\theta-2}\|u\|_{\bs+2,\mathcal{T}_h} + h^{\min\{\theta-2,2\}}\left[h^{\mu-1} \|u\|_{s,\mathcal{T}_h} +
      h^{\bar{\mu}-1}
      \|u\|_{\bs+2,\mathcal{T}_h}\right]\right).\label{eq:error8}
    \end{align}
    Combining estimates~\eqref{eq:error7} and \eqref{eq:error8} and using the triangle inequality, we obtain the desired result for $|||e_u|||$.
  \end{proof}

\end{theorem}
Finally using Theorem~\ref{theorem2} and Theorem~\ref{theorem3}, we
can prove the following estimate for $\|e_u\|$.
\begin{theorem}
\label{theorem4}
  Let $u_h,v_h$ satisfy \eqref{eq:8}. There exists some positive
  constant $C$ such that
  \begin{align*}
    \|e_u\| \le Ch^{\min\{\mu+\theta-3,\mu\}}\|u\|_{s,\mathcal{T}_h}
    + Ch^{\min\{\bar{\mu}+\theta-3,\bar{\mu}\}}\|u\|_{\bs+2,\mathcal{T}_h},
  \end{align*}
  provided the conditions on $\sigma_0$ in Theorem~\ref{theorem1} holds.
\end{theorem}

\section{Numerical Experiments}
\label{sec:3}
In this section, we present some numerical experiments to validate the
theoretical results. To perform the grid refinement analysis, we
choose the exact solution $u(x,y) = 1000\,x^4y^4(1-x)^4(1-y)^4$ and
calculate the right hand side. We then calculate the DGFEM solution $u_h(x,y)$
and compute the error in the $L^2$ and the energy norms. The numerical
experiments were conducted using FreeFem++. We consider two cases
with different choices of penalty parameters and compare the results of the
numerical scheme. In all the numerical results, we set the parameter
$\sigma_0 = 1$.\\

\begin{table}
  \centering
  \begin{tabular}{|p{1cm}|p{2cm}|p{2cm}|p{2cm}|p{2.5cm}|p{2.5cm}|}
    \hline
    $p$ & $h$ & $\|u-u_h\|$ & $|||u-u_h|||$ & O($\|u-u_h\|$) &
                                                               O($|||u-u_h|||$)\\
    \hline
    1 & 0.141421 & 0.000419 & 0.011011 & -  & - \\
        & 0.094281 & 0.000375 & 0.007711 & 0.270249 & 0.878539 \\
        & 0.070711 & 0.000371 & 0.006441 & 0.036169 & 0.625625 \\
        & 0.056569 & 0.000372 & 0.005790 & -0.009047 & 0.477075 \\
        & 0.047140 & 0.000373 & 0.005412 & -0.017329 & 0.371163 \\
        & 0.040406 & 0.000374 & 0.005172 & -0.017346 & 0.293645 \\
    \hline
    2 & 0.282843 & 0.000402 & 0.017695 & -  & - \\
        & 0.141421 & 0.000027 & 0.002197 & 3.887976 & 3.009653 \\
        & 0.094281 & 0.000007 & 0.000960 & 3.183386 & 2.043090 \\
        & 0.070711 & 0.000003 & 0.000544 & 2.876605 & 1.973371 \\
        & 0.056569 & 0.000002 & 0.000350 & 2.695863 & 1.977692 \\
        & 0.047140 & 0.000001 & 0.000244 & 2.565206 & 1.983976 \\
    \hline
    3 & 0.282843 & 0.000034 & 0.002809 & -  & - \\
        & 0.141421 & 0.000002 & 0.000240 & 4.168647 & 3.550827 \\
        & 0.094281 & 0.000000 & 0.000068 & 3.915067 & 3.120461 \\
        & 0.070711 & 0.000000 & 0.000028 & 3.929946 & 3.049651 \\
        & 0.056569 & 0.000000 & 0.000014 & 3.947376 & 3.031438 \\
        & 0.047140 & 0.000000 & 0.000008 & 3.960072 & 3.024891 \\
    \hline
  \end{tabular}
  \caption{Error values and rates of convergence of the Mixed DGFEM for the biharmonic equation. The penalty
    parameter is chosen to be $\alpha_k = \sigma_0|e_k|^{-3}p^2$ with
    $\sigma_0 = 1$. The
    rate of convergence for the linear DGFEM is not significant and a
    sub-optimal convergence for the piecewise quadratic elements is
    observed which is an observation made by \cite{gudi2008mixed}.}
  \label{tab:1}
\end{table}
First, we compute the rates of convergence for the DGFEM considered by \cite{gudi2008mixed}. From Table~\ref{tab:1} we observe that the rates of
convergence for the linear DGFEM is not significant and a sub-optimal
convergence (rate of convergence decreases rapidly) in the $L^2$--norm
was observed for piecewise quadratic
elements. The penalty parameter was chosen as $\alpha_k = \sigma_0|e_k|^{-3}p^2$
for the mixed DGFEM. However, the choice of the penalty parameter is
crucial for the mixed method as choosing a higher value may result in
sub-optimal convergence.\\

\begin{figure}
  \begin{subfigure}{0.5\textwidth}
    \includegraphics[width=\textwidth]{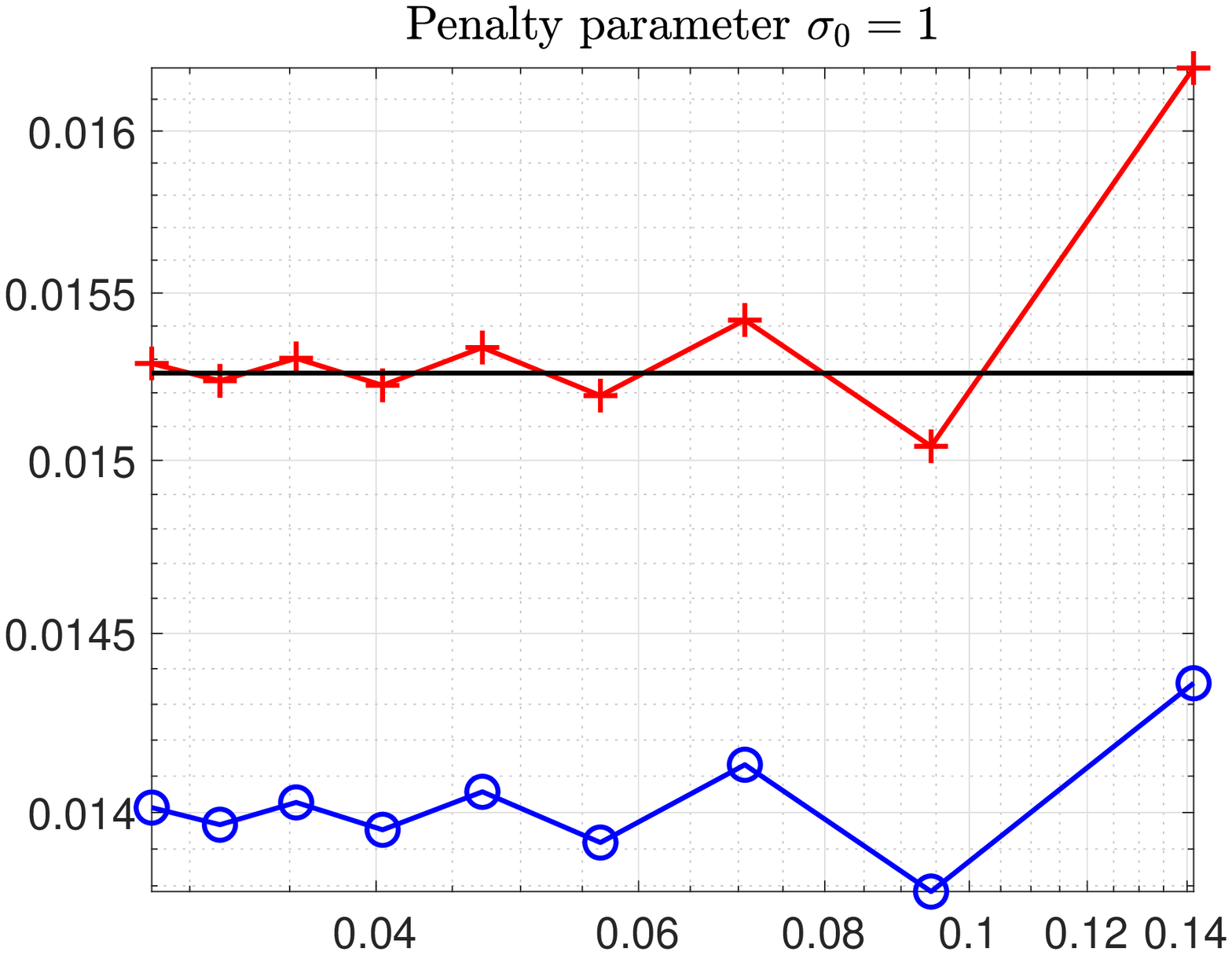}
    \caption{$\sigma_0=1$}
  \end{subfigure}
  \begin{subfigure}{0.5\textwidth}
    \includegraphics[width=\textwidth]{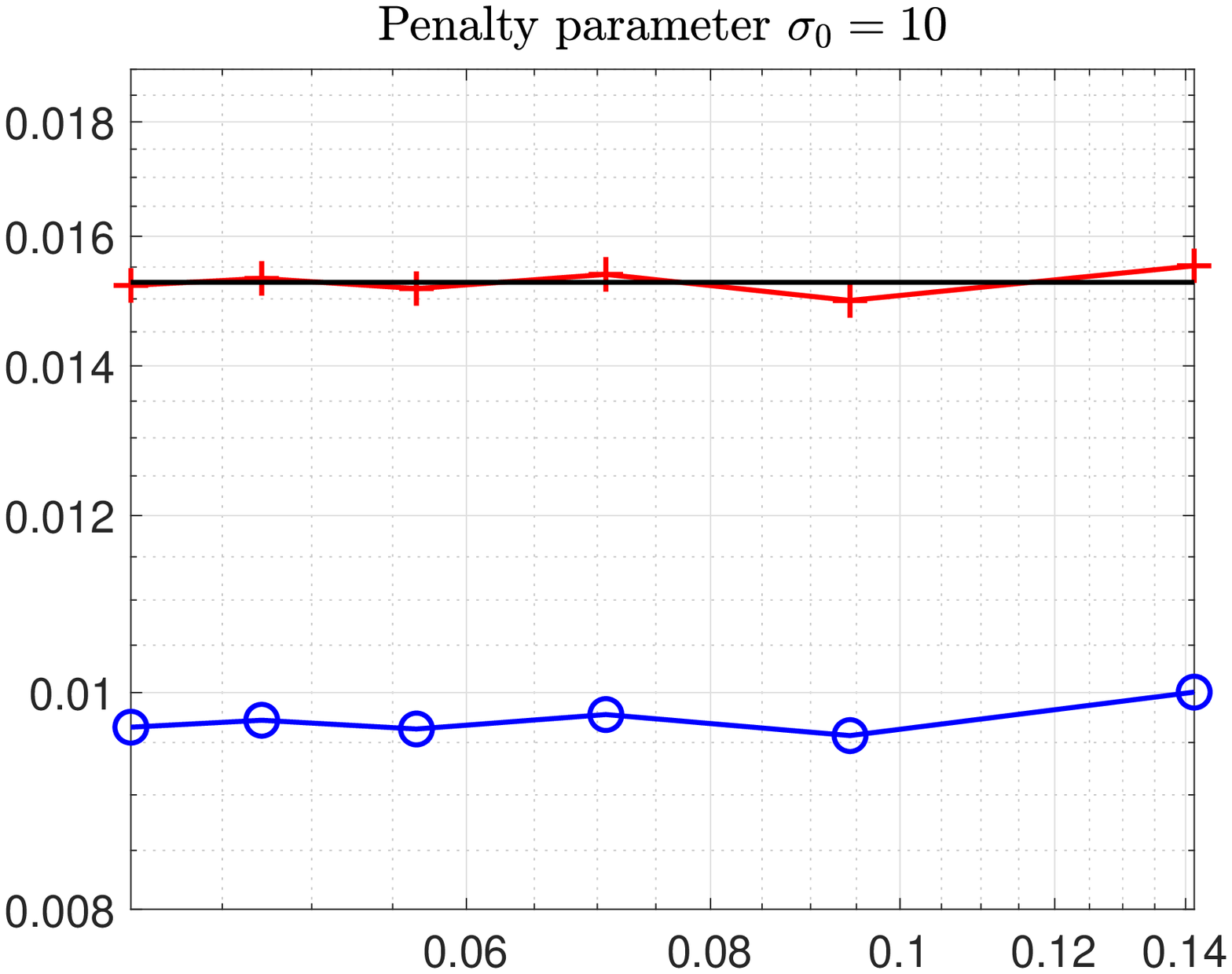}
    \caption{$\sigma_0=10$}
  \end{subfigure}
  \caption{Convergence of the solution in terms of the maximum value
    of the approximate solution ($y$-axis) with the mesh parameter
    ($x$-axis) for piecewise linear polynomials. The red curve (+) denotes maximum value
    of the DG solution with $\alpha_k = \sigma_0 |e_k|^{-1}$ and the
    blue curve (O) denotes the maximum value of the DG solution with
    $\alpha_k = \sigma_0 |e_k|^{-3}$ considered by
    \cite{gudi2008mixed}. The horizontal black line denotes the
    maximum value of the exact solution. We observe that the
    convergence of the solution is observed only in the red curve
    (lower penalty term), whereas convergence is not observed in the
    blue curve (higher penalty term).}
  \label{fig:1}
\end{figure}

\begin{table}
  \centering
  \begin{tabular}{|p{1cm}|p{2cm}|p{2cm}|p{2cm}|p{2.5cm}|p{2.5cm}|}
    \hline
    $p$ & $h$ & $\|u - u_h\|$ & $|||u-u_h|||$ & O($\|u-u_h\|$) &
                                                                 O($|||u-u_h|||$)\\
    \hline
    1   & 0.141421 & 0.004153 & 0.135755 & -  & - \\
        & 0.094281 & 0.001258 & 0.060027 & 2.944717 & 2.012619 \\
        & 0.070711 & 0.000549 & 0.034655 & 2.882794 & 1.909588 \\
        & 0.056569 & 0.000291 & 0.022947 & 2.839182 & 1.847442 \\
        & 0.047140 & 0.000174 & 0.016491 & 2.814693 & 1.811970 \\
        & 0.040406 & 0.000113 & 0.012517 & 2.799510 & 1.788672 \\
    \hline
    2 & 0.282843 & 0.006098 & 0.361025 & -  & - \\
        & 0.141421 & 0.000578 & 0.071355 & 3.400331 & 2.339021 \\
        & 0.094281 & 0.000142 & 0.026470 & 3.452994 & 2.445662 \\
        & 0.070711 & 0.000053 & 0.013065 & 3.462819 & 2.454307 \\
        & 0.056569 & 0.000024 & 0.007532 & 3.473787 & 2.468250 \\
        & 0.047140 & 0.000013 & 0.004798 & 3.478131 & 2.473979 \\
    \hline
    3  & 0.282843 & 0.000159 & 0.018969 & -  & - \\
        & 0.141421 & 0.000013 & 0.003115 & 3.562213 & 2.606456 \\
        & 0.094281 & 0.000003 & 0.001173 & 3.409106 & 2.409145 \\
        & 0.070711 & 0.000001 & 0.000522 & 3.817537 & 2.811388 \\
        & 0.056569 & 0.000000 & 0.000267 & 4.012594 & 3.007459 \\
        & 0.047140 & 0.000000 & 0.000151 & 4.124334 & 3.121269 \\
    \hline
  \end{tabular}
  \caption{Error values and rates of convergence for the mixed DGFEM
    for the biharmonic equation. The rate of convergence is close to
    $p+1$ in the $L^2$-norm and close to $p$ in the energy norm
    indicating optimal rates of convergence. The
    penalty parameter is chosen to be $\alpha_k = \sigma_0
    |e_k|^{-1}p^2$ with $\sigma_0 = 1$.}
  \label{tab:2}
\end{table}
In Table~\ref{tab:2}, we present the
convergence results by choosing a lower value of the penalty parameter
with $\alpha_k = \sigma_0|e_k|^{-1}p^2$. The constant $\sigma_0$ in
all the cases was set to be equal to $1$.
\begin{figure}
  \centering
  \includegraphics[width=\textwidth,height=0.4\textheight]{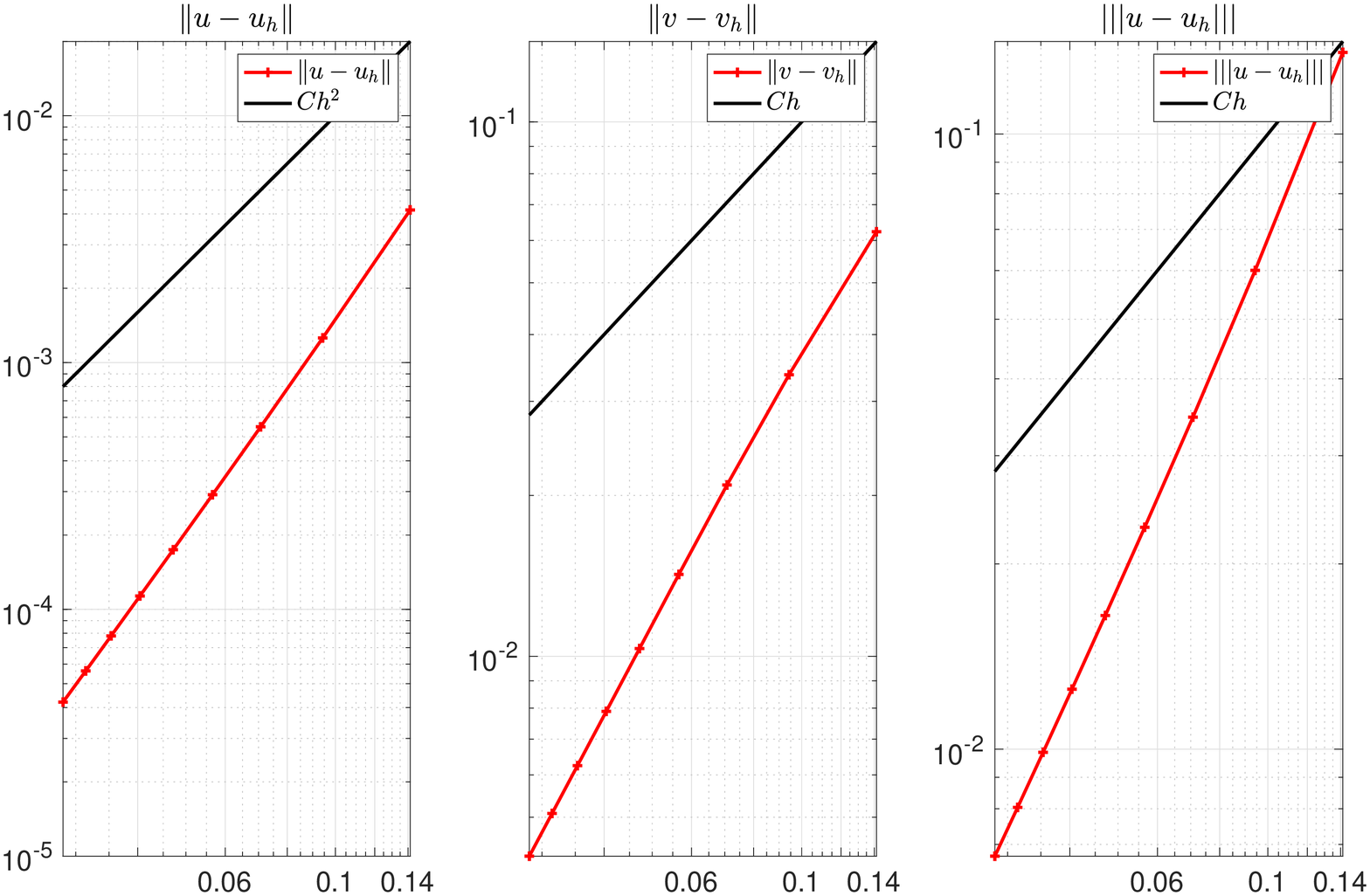}
  \caption{Rates of convergence of various errors ($y$-axes) with mesh
    size ($x$-axes) for piecewise linear Discontinuous
    Galerkin Method. The black solid line shows the expected rate of convergence predicted by the theoretical results.}
  \label{fig:2}
\end{figure}
We observe a significant
improvement in the convergence rates for piecewise linear elements and
for piecewise quadratic elements. We observe
that on choosing a lesser value of penalty, the solution becomes more
accurate with refinement. This is especially strong in the linear
case, where the solution converges rapidly (rate close to $\approx
2.5$) to the exact solution and the magnitude of the $L^2-$error
$\|u-u_h\|$ is significantly lower in the last iteration. This is
illustrated in Figure~\ref{fig:1} where convergence
is observed only when $\alpha_k = \sigma_0|e_k|^{-1}$. However, the
theoretical results in Theorem~\ref{theorem3}~and~\ref{theorem4} predict optimal estimates for the energy norm
$|||e_u|||$ $(\approx o(1))$ and sub-optimal estimates for $\|e_u\|$
$(\approx o(1))$, repectively, despite the observed higher convergence rates in the case of $\|e_u\|$.
\begin{figure}
  \includegraphics[width=\textwidth,height=0.4\textheight]{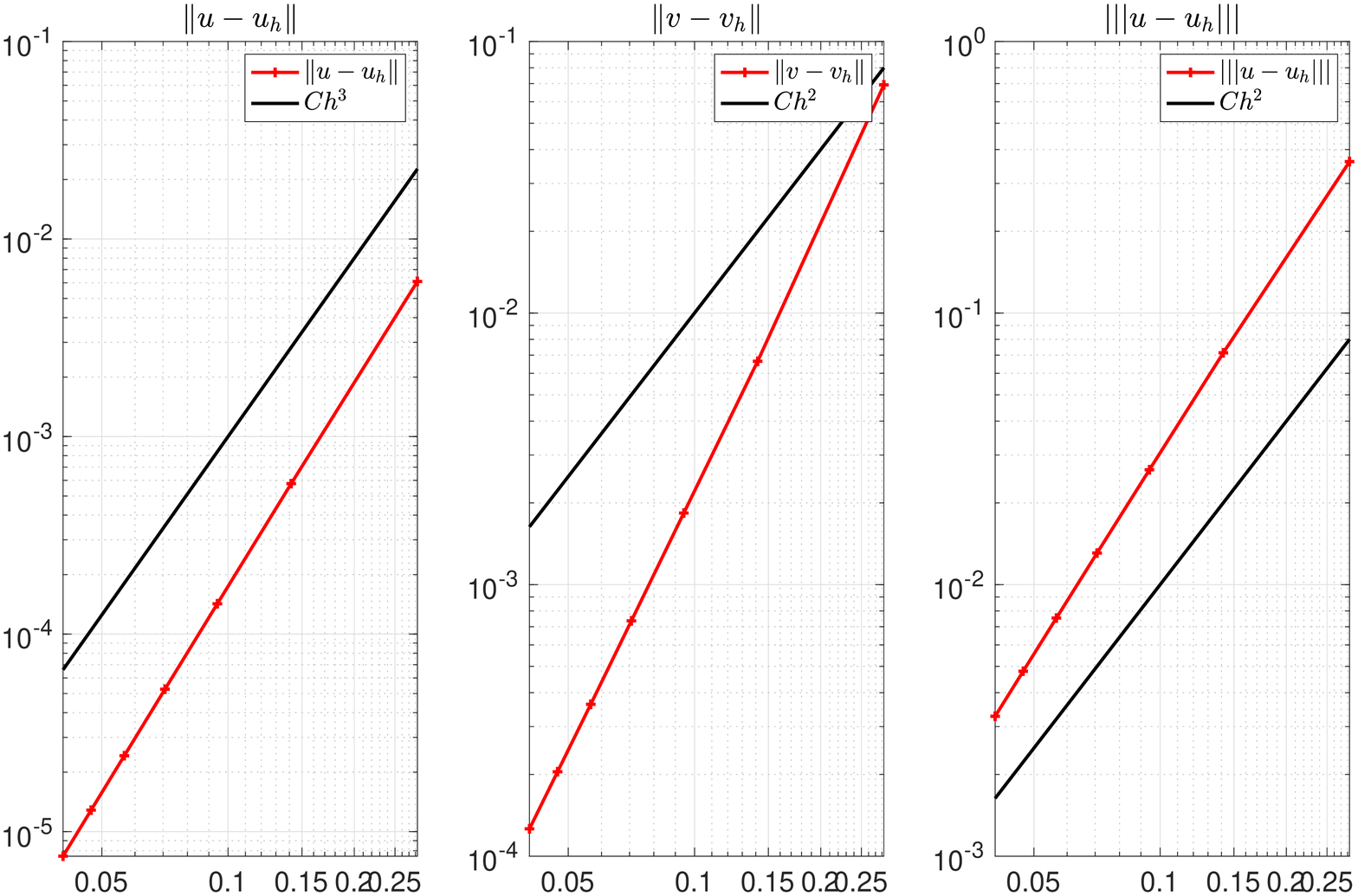}
  \caption{Same as Figure~\ref{fig:2} but for piecewise quadratic Discontinuous
    Galerkin Method.}
  \label{fig:3}
\end{figure}
We observe a higher convergence rate for
piecewise quadratic elements in the $L^2$ and energy norms while the
error magnitude is comparable to the previous case in
Table~\ref{tab:1}. The higher convergence rates are also predicted by
the theoretical estimates for piecewise quadratic case. We observe
that optimal convergence rates are preserved for the piecewise cubic
case, although the error magnitudes are higher than that observed from
Table~\ref{tab:1}. The theoretical error estimates predict optimal convergence rates for the errors $|||e_u|||$ and $\|e_u\|$ for both piecewise quadratic and cubic cases.\\
\begin{figure}
  \includegraphics[width=\textwidth,height=0.4\textheight]{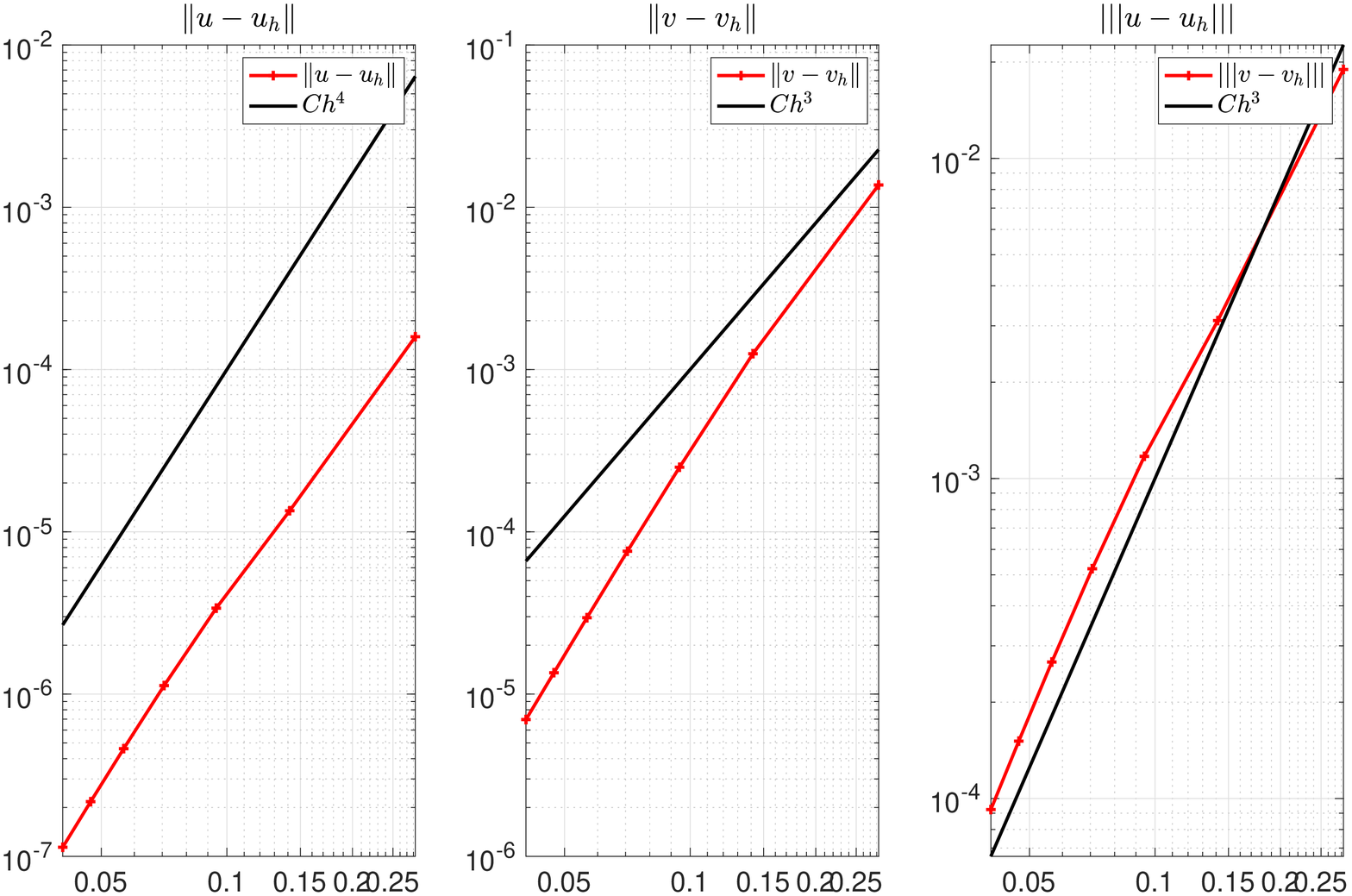}
  \caption{Same as Figure~\ref{fig:2} but for piecewise cubic Discontinuous
    Galerkin Method.}
  \label{fig:4}
\end{figure}

Similar observations were made for the convergence in the auxiliary
variable $v$ which is summarized in Table~\ref{tab:3}. \cite{gudi2008mixed} observed a sub--optimal convergence rate
($\approx p-1$) for the auxiliary variable in the $L^2$ norm. This was
not observed in the case when the penalty parameter was chosen as
$\alpha_k = \sigma_0|e_k|^{-1}p^2$ and the convergence rates are
optimal $(\approx p+1)$ in the $L^2$ norm. The observed rate of convergence is better
than the theoretically established result in Theorem~\ref{theorem1},
which predicts a sub--optimal convergence in $v$ $(\approx p)$. We observe that the
current choice of penalty term works well to approximate the auxiliary
variable even for the piecewise linear case. Figures~\ref{fig:2},~\ref{fig:3}~and~\ref{fig:4} shows the decay in the approximation error with respect to the mesh parameter $h$ for piecewise linear, quadratic and cubic polynomial approximations, respectively. 
\begin{table}
  \centering
  \begin{tabular}{|p{2cm}|p{2cm}|p{2cm}|p{2cm}|p{2cm}|p{2cm}|}
    \hline
    \multicolumn{2}{|c|}{} & \multicolumn{2}{|c|}{$\alpha_k = \sigma_0
                             |e_k|^{-1}p^2$} & \multicolumn{2}{|c|}{$\alpha_k
                                               = \sigma_0
                                               |e_k|^{-3}p^2$
                                               }\\
    \hline
    $p$ & $N$ & $\| v - v_h\|$ & O($\| v - v_h \|$) & $\| v - v_h\|$ &
                                                                       O($\| v - v_h\|$)\\
    \hline
    1 & 0.141421 & 0.062298 & - & 0.030978 & - \\
                           & 0.094281 & 0.033646 & 1.519314 & 0.022990
                               & 0.735562\\
                           & 0.070711 & 0.020921 & 1.651590 & 0.024028
                               & -0.153535 \\
                           & 0.056569 & 0.014233 & 1.726174 & 0.025236
                               & -0.219860 \\
                           & 0.047140 & 0.010340 & 1.752746 & 0.026058
                               & -0.175885 \\
                           & 0.040406 & 0.007888 & 1.755539 & 0.026604
                               & -0.134564 \\
    \hline
    2 & 0.282843 & 0.069201 & - & 0.038757 & -\\
                           & 0.141421 & 0.006641 & 3.381222 & 0.005352
                               & 2.856424\\
                           & 0.094281 & 0.001834 & 3.173283 & 0.002105
                               & 2.301325\\
                           & 0.070711 & 0.000736 & 3.176304 & 0.001346
                               & 1.555525\\
                           & 0.056569 & 0.000362 & 3.171694 & 0.001035
                               & 1.176948\\
                           & 0.047140 & 0.000205 & 3.136815 & 0.000855
                               & 1.044048\\
    \hline
    3 & 0.282843 & 0.013709 & - & 0.006130 & -\\
                           & 0.141421 & 0.001251 & 3.454083 & 0.000502
                               & 3.609629\\
                           & 0.094281 & 0.000250 & 3.972908 & 0.000288
                               & 1.370246\\
                           & 0.070711 & 0.000076 & 4.143599 & 0.000188
                               & 1.479595\\
                           & 0.056569 & 0.000030 & 4.227323 & 0.000128
                               & 1.722811\\
                           & 0.047140 & 0.000014 & 4.278536 & 0.000091
                               & 1.869907\\
    \hline
  \end{tabular}
  \caption{The rates of convergence for the auxiliary variable for the
    DGFEM with respect to the mesh size and penalty parameter. We
    observe optimal convergence rate for the auxiliary variable in the
    $L^2$ and the energy norm for the current choice of penalty
    parameter and suboptimal rate of $p-1$ for the choice in
    \cite{gudi2008mixed}. We set $\sigma_0 = 1.$}
  \label{tab:3}
\end{table}

\section{Conclusion}
In this work, we have considered a mixed Discontinuous Galerkin Finite
Element Method to solve the biharmonic equation subject to clamped
boundary conditions. We have derived the weak formulation of the
problem and established error estimates for $h$-refinement. The
theoretical results predict optimal convergence rates in the energy
norm for the primal variable $u$, whereas optimal convergence is
predicted for the piecewise quadratic and cubic case in the
$L^2$-norm. We performed a
series of numerical experiments using FreeFem++, and verified the
theoretical results. We observed that the choice of the penalty term
is crucial and must be chosen to be of the form $\alpha_k = \sigma_0
|e_k|^{-1}p^2$ to obtain optimal error estimates in
$h$-refinement. Significant improvements in convergence rates for the
piecewise linear and quadratic elements were observed as a result.

\bibliographystyle{IMANUM-BIB}
\bibliography{biblio}

\begin{thebibliography}{}

\bibitem[Arnold {\em et~al.}(2002)Arnold, Brezzi, Cockburn, \&
  Marini]{arnold2002unified}
{\sc Arnold, D.~N., Brezzi, F., Cockburn, B. \& Marini, L.~D.} (2002)
\newblock {Unified Analysis of Discontinuous Galerkin Methods for Elliptic
  Problems}.
\newblock {\em SIAM Journal on Numerical Analysis\/}, {\bf 39}, 1749--1779.

\bibitem[Brenner \& Scott(2007)Brenner \& Scott]{brenner2007mathematical}
{\sc Brenner, S.~C. \& Scott, R.} (2007)
\newblock {\em {The Mathematical Theory of Finite Element Methods}\/},
  vol.~15.
\newblock Springer Science \& Business Media.

\bibitem[Brenner \& Sung(2005)Brenner \& Sung]{brenner2005c}
{\sc Brenner, S.~C. \& Sung, L.-Y.} (2005)
\newblock {$C^0$ Interior Penalty Methods for Fourth Order Elliptic Boundary
  Value Problems on Polygonal Domains}.
\newblock {\em Journal of Scientific Computing\/}, {\bf 22}, 83--118.

\bibitem[Ciarlet \& Raviart(1974)Ciarlet \& Raviart]{ciarlet1974mixed}
{\sc Ciarlet, P.~G. \& Raviart, P.-A.} (1974)
\newblock {A Mixed Finite Element Method for the Biharmonic Equation}.
\newblock {\em Mathematical Aspects of Finite Elements in Partial Differential
  Equations\/}.
\newblock Elsevier, pp. 125--145.

\bibitem[Cockburn {\em et~al.}(2012)Cockburn, Karniadakis, \&
  Shu]{cockburn2012discontinuous}
{\sc Cockburn, B., Karniadakis, G.~E. \& Shu, C.-W.} (2012)
\newblock {\em {Discontinuous Galerkin Methods: Theory, Computation and
  Applications}\/},  vol.~11.
\newblock Springer Science \& Business Media.

\bibitem[Danumjaya \& Pani(2012)Danumjaya \& Pani]{danumjaya2012mixed}
{\sc Danumjaya, P. \& Pani, A.~K.} (2012)
\newblock {Mixed Finite Element Methods for a Fourth Order Reaction Diffusion
  Equation}.
\newblock {\em Numerical Methods for Partial Differential Equations\/}, {\bf
  28}, 1227--1251.

\bibitem[Engel {\em et~al.}(2002)Engel, Garikipati, Hughes, Larson, Mazzei, \&
  Taylor]{engel2002continuous}
{\sc Engel, G., Garikipati, K., Hughes, T., Larson, M., Mazzei, L. \& Taylor,
  R.~L.} (2002)
\newblock {Continuous/Discontinuous Finite Element Approximations of
  Fourth-Order Elliptic Problems in Structural and Continuum Mechanics with
  Applications to Thin Beams and Plates, and Strain Gradient Elasticity}.
\newblock {\em Computer Methods in Applied Mechanics and Engineering\/}, {\bf
  191}, 3669--3750.

\bibitem[Gudi {\em et~al.}(2008)Gudi, Nataraj, \& Pani]{gudi2008mixed}
{\sc Gudi, T., Nataraj, N. \& Pani, A.~K.} (2008)
\newblock {Mixed Discontinuous Galerkin Finite Element Method for the
  Biharmonic Equation}.
\newblock {\em Journal of Scientific Computing\/}, {\bf 37}, 139--161.

\bibitem[Monk(1987)Monk]{monk1987mixed}
{\sc Monk, P.} (1987)
\newblock A mixed finite element method for the biharmonic equation.
\newblock {\em SIAM Journal on Numerical Analysis\/}, {\bf 24}, 737--749.

\bibitem[Rivi\'ere(2008)Rivi\'ere]{riviere2008discontinuous}
{\sc Rivi\'ere, B.} (2008)
\newblock {\em {Discontinuous Galerkin Methods for Solving Elliptic and
  Parabolic Equations: Theory and Implementation}\/}.
\newblock SIAM.

\end{thebibliography}

\end{document}